\documentclass{article}
\usepackage[margin=1.25in]{geometry}
\usepackage{amsmath}
\usepackage{amssymb}
\usepackage{graphicx}
\usepackage{cite}
\usepackage{float}

\newcommand{\define}[1]{\emph{#1}}

\newcommand{\citeasnoun}[1]{Ref.~\citenum{#1}}

\newcommand{\PUNT}[1]{}

\usepackage{algorithmic}
\newfloat{algorithm}{tbp}{loa}
\floatname{algorithm}{Algorithm}

\newcommand{\secref}[1]{section~\ref{sec:#1}}
\newcommand{\secreftwo}[2]{sections \ref{sec:#1} and~\ref{sec:#2}}

\newcommand{\figref}[1]{fig.~\ref{fig:#1}}

\newcommand{\algref}[1]{algorithm~\ref{alg:#1}}

\newcommand{\eqnumref}[1]{(\ref{eq:#1})}

\renewcommand{\eqref}[1]{eq.~\eqnumref{#1}}

\newcommand{\onedim}{1d}

\newcommand{\DFT}{DFT}
\newcommand{\DFTs}{DFTs}

\newcommand{\DCT}{DCT}

\newcommand{\DST}{DST}

\newcommand{\FFT}{FFT}
\newcommand{\FFTs}{FFTs}
\newcommand{\DIT}{DIT}
\newcommand{\DIF}{DIF}
\newcommand{\FFTW}{FFTW}
\newcommand{\FFTWthree}{FFTW}

\newcommand{\SPIRAL}{SPIRAL}

\newcommand{\union}{\cup}

\newcommand{\imagunit}{i}

\newcommand{\rootunity}{\omega}

\newcommand{\IODIM}{I/O dimension}
\newcommand{\IODIMs}{{\IODIM}s}
\newcommand{\iodim}[3]{(#1,#2,#3)}
\newcommand{\aiodim}{d}
\newcommand{\n}{n}
\newcommand{\is}{\iota}
\newcommand{\os}{o}

\newcommand{\IOTENS}{I/O tensor}
\newcommand{\IOTENSs}{{\IOTENS}s}
\newcommand{\aiotens}{t}
\newcommand{\iotens}[1]{\left\{#1\right\}}
\newcommand{\iotensone}{\iotens{}}

\newcommand{\dftprob}[4]{\operatorname{dft}(#1, #2, #3, #4)}

\newcommand{\dftx}{{\mathbf X}}
\newcommand{\dfty}{{\mathbf Y}}
\newcommand{\dftn}{{\mathbf N}}
\newcommand{\dftv}{{\mathbf V}}
\newcommand{\dfti}{{\mathbf I}}
\newcommand{\dfto}{{\mathbf O}}

\newcommand{\rank}{\rho}
\newcommand{\radix}{r}
\newcommand{\m}{m}
\newcommand{\rankof}[1]{\left| #1 \right|}
\newcommand{\ndx}{k}  

\newcommand{\copyo}[1]{\operatorname{copy-o}(#1)}
\newcommand{\suchthat}{\mid}
\newcommand{\genfft}{\texttt{genfft}}

\usepackage{fancyhdr}
\newcommand{\BookTitle}{\textit{Fast Fourier Transforms} (C.S.~Burrus, ed.), ch.~11, Rice Univ.: Connexions, 2008.}

\fancypagestyle{plain}{
  \fancyhf{}
  \fancyhead[C]{\small \BookTitle}
  \fancyfoot[C]{\thepage}
  
}

\usepackage{authblk}


\begin{document}

\title{Implementing {\FFTs} in Practice}
\author[1]{Steven G. Johnson}
\author[2]{Matteo Frigo}
\affil[1]{Department of Mathematics, Massachusetts Institute of Technology.}
\affil[2]{Cilk Arts, Inc.}

\date{Published 2008 as ch.~11 in \textit{Fast Fourier Transforms}, ed.~C.~S.~Burrus (Rice). \\[0.5em]
Produced by The Connexions Project \\ 
and licensed under the Creative Commons Attribution License. }
\maketitle

\section{Introduction}

Although there are a wide range of fast Fourier transform ({\FFT})
algorithms, involving a wealth of mathematics from number theory to
polynomial algebras, the vast majority of {\FFT} implementations in
practice employ some variation on the Cooley--Tukey
algorithm~\cite{CooleyTu65}.  The Cooley--Tukey algorithm can be
derived in two or three lines of elementary algebra.  It can be
implemented almost as easily, especially if only power-of-two sizes
are desired; numerous popular textbooks list short {\FFT} subroutines
for power-of-two sizes, written in the language {\it du jour}.  The
implementation of the Cooley--Tukey algorithm, at least, would
therefore seem to be a long-solved problem.  In this chapter, however,
we will argue that matters are not as straightforward as they might
appear.

For many years, the primary route to improving upon the Cooley--Tukey
{\FFT} seemed to be reductions in the count of arithmetic operations,
which often dominated the execution time prior to the ubiquity of fast
floating-point hardware.  Therefore, great effort was expended towards
finding new algorithms with reduced arithmetic
counts~\cite{DuhamelVe90}, from Winograd's method to achieve
$\Theta(n)$ multiplications\footnote{We employ the standard asymptotic
notation of $O$ for asymptotic upper bounds, $\Theta$ for asymptotic
tight bounds, and $\Omega$ for asymptotic lower
bounds~\cite{TAOCP-I}.} (at the cost of many more
additions)~\cite{Winograd78,Heideman86,Duhamel90,DuhamelVe90} to the
split-radix variant on Cooley--Tukey that long achieved the lowest
known total count of additions and multiplications for power-of-two
sizes~\cite{Yavne68,Duhamel84,Vetterli84,Martens84,DuhamelVe90} (but
was recently improved upon~\cite{Johnson07,Lundy07}).  The question of
the minimum possible arithmetic count continues to be of fundamental
theoretical interest---it is not even known whether better than
$\Theta(n \log n)$ complexity is possible, since $\Omega(n \log n)$
lower bounds on the count of additions have only been proven subject
to restrictive assumptions about the
algorithms~\cite{Morgenstern73,Pan86,Papadimitriou79}.  Nevertheless,
the difference in the number of arithmetic operations, for
power-of-two sizes $n$, between the 1965 radix-2 Cooley--Tukey
algorithm ($\sim 5 n \log_2 n$~\cite{CooleyTu65}) and the currently
lowest-known arithmetic count ($\sim \frac{34}{9} n \log_2
n$~\cite{Johnson07,Lundy07}) remains only about 25\%.

\begin{figure}
\centering\includegraphics[width=0.90\columnwidth]{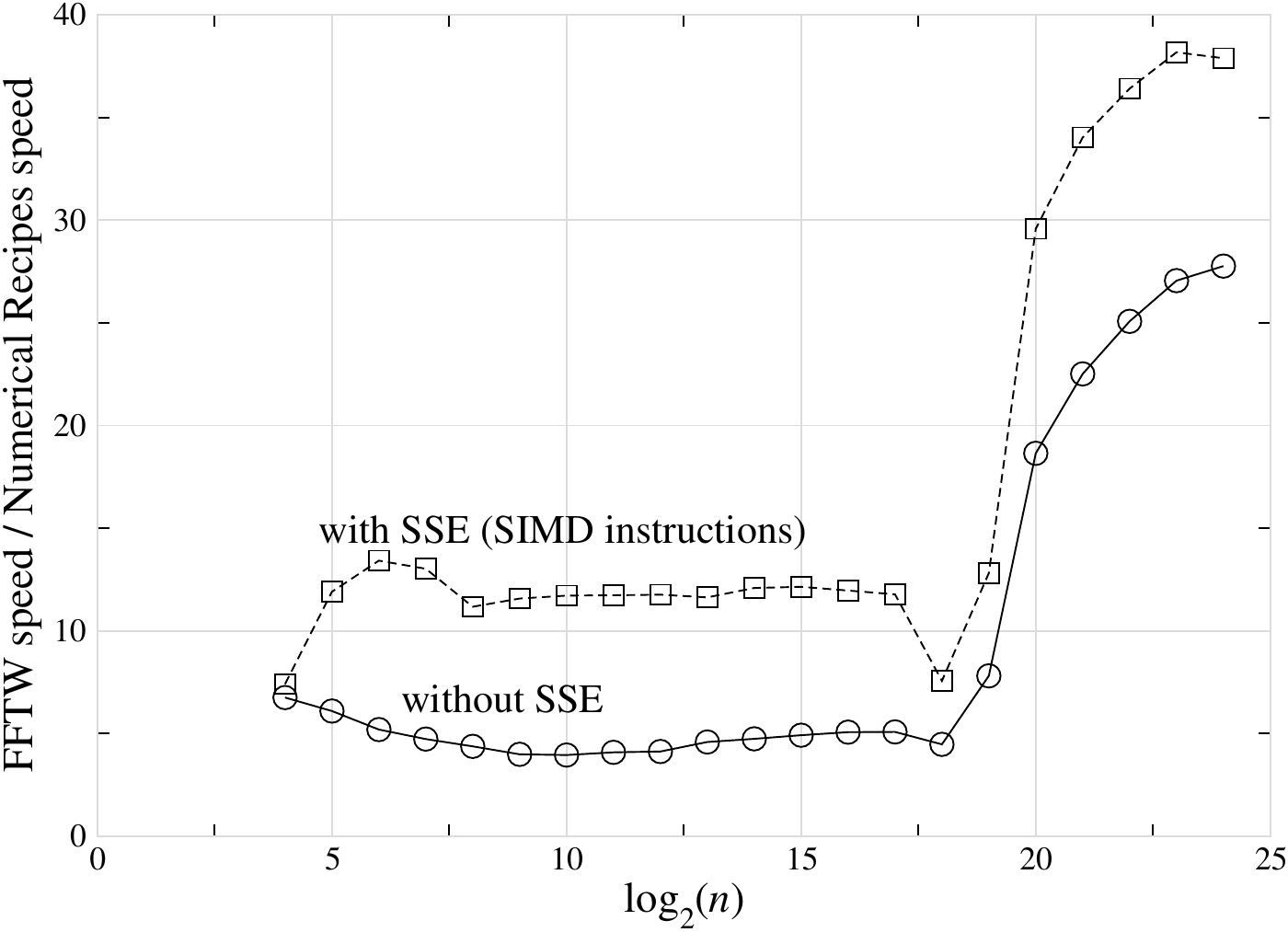}
\caption{The ratio of speed (1/time) between a highly optimized FFT
({\FFTW} 3.1.2~\cite{FFTWweb,FFTW05}) and a typical textbook radix-2
implementation ({\it Numerical Recipes in C}~\cite{PressFlaTeu92}) on a 3~GHz Intel
Core Duo with the Intel C compiler 9.1.043, for single-precision
complex-data {\DFTs} of size $n$, plotted versus $\log_2 n$.  Top line
(squares) shows {\FFTW} with SSE SIMD instructions enabled, which
perform multiple arithmetic operations at once (see \secref{genfft:simd});
bottom line (circles) shows {\FFTW} with SSE disabled, which thus
requires a similar number of arithmetic instructions to the textbook
code.  (This is not intended as a criticism of {\it Numerical
Recipes}---simple radix-2 implementations are reasonable for
pedagogy---but it illustrates the radical differences between
straightforward and optimized implementations of {\FFT} algorithms,
even with similar arithmetic costs.)  For $n \gtrsim 2^{19}$, the
ratio increases because the textbook code becomes much slower (this
happens when the {\DFT} size exceeds the level-2 cache).}
\label{fig:fftwnr}
\end{figure}

And yet there is a vast gap between this basic mathematical theory and
the actual practice---highly optimized {\FFT} packages are often an
order of magnitude faster than the textbook subroutines, and the
internal structure to achieve this performance is radically different
from the typical textbook presentation of the ``same'' Cooley--Tukey
algorithm.  For example, \figref{fftwnr} plots the ratio of benchmark
speeds between a highly optimized {\FFT}~\cite{FFTWweb,FFTW05} and a
typical textbook radix-2 implementation~\cite{PressFlaTeu92}, and the
former is faster by a factor of 5--40 (with a larger ratio as $n$
grows).  Here, we will consider some of the reasons for this
discrepancy, and some techniques that can be used to address the
difficulties faced by a practical high-performance {\FFT}
implementation.\footnote{We won't address the question of
parallelization on multi-processor machines, which adds even greater
difficulty to {\FFT} implementation---although multi-processors are
increasingly important, achieving good serial performance is a basic
prerequisite for optimized parallel code, and is already hard enough!}

In particular, in this chapter we will discuss some of the lessons
learned and the strategies adopted in the {\FFTW} library.  {\FFTW}~\cite{FFTWweb,FFTW05} is a widely used free-software library that computes the
discrete Fourier transform ({\DFT}) and its various special cases.
Its performance is competitive even with manufacturer-optimized programs~\cite{FFTW05},
and this performance is \emph{portable} thanks the structure of the
algorithms employed, self-optimization techniques, and highly
optimized kernels ({\FFTW}'s \define{codelets}) generated by a
special-purpose compiler.

This chapter is structured as follows.  First (\secref{fftintro}), we
briefly review the basic ideas behind the Cooley--Tukey algorithm and
define some common terminology, especially focusing on the many
degrees of freedom that the abstract algorithm allows to
implementations.  Second (\secref{cache}), we consider a basic
theoretical model of the computer memory hierarchy and its impact on
{\FFT} algorithm choices: quite general considerations push
implementations towards large radices and explicitly recursive
structure.  Unfortunately, general considerations are not sufficient
in themselves, so we will explain in \secref{struct} how {\FFTW}
self-optimizes for particular machines by selecting its algorithm at
runtime from a composition of simple algorithmic steps.  Furthermore,
\secref{genfft} describes the utility and the principles of automatic
code generation used to produce the highly optimized building blocks
of this composition, {\FFTW}'s codelets.  Finally, we will
consider some non-performance issues involved in a practical {\FFT}
implementation, and in particular the question of accuracy
(\secref{accuracy}) and the crucial consideration of generality
(\secref{generality}).

\section{Review of the Cooley--Tukey {\FFT}}
\label{sec:fftintro}

The (forward, one-dimensional)
discrete Fourier transform ({\DFT}) of an array $\dftx$ of $n$ complex numbers is
the array $\dfty$ given by
\begin{equation}
\dfty[k] = \sum_{\ell = 0}^{n-1} \dftx[\ell] \rootunity_n^{\ell k} \ ,
\label{eq:dft}
\end{equation}
where $0 \leq k < n$ and $\rootunity_n = \exp(-2\pi \imagunit /n)$ is a
primitive root of unity.  Implemented directly, \eqref{dft} would
require $\Theta(n^2)$ operations; fast Fourier transforms are $O(n
\log n)$ algorithms to compute the same result.  The most important
{\FFT} (and the one primarily used in {\FFTW}) is known as the
``Cooley--Tukey'' algorithm, after the two authors who rediscovered and
popularized it in 1965~\cite{CooleyTu65}, although it had been
previously known as early as 1805 by Gauss as well as by later
re-inventors~\cite{Heideman84}.  The basic idea behind this {\FFT} is
that a {\DFT} of a composite size $n = n_1 n_2$ can be re-expressed in
terms of smaller {\DFTs} of sizes $n_1$ and $n_2$---essentially, as a
two-dimensional {\DFT} of size $n_1 \times n_2$ where the output is
\emph{transposed}.  The choices of factorizations of $n$, combined
with the many different ways to implement the data re-orderings of the
transpositions, have led to numerous implementation strategies for the
Cooley--Tukey {\FFT}, with many variants distinguished by their own
names~\cite{DuhamelVe90,VanLoan92}.  {\FFTW} implements a space of
\emph{many} such variants, as described in \secref{struct}, but here
we derive the basic algorithm, identify its key features, and outline
some important historical variations and their relation to {\FFTW}.

The Cooley--Tukey algorithm can be derived as follows.  If $n$ can be
factored into $n = n_1 n_2$, \eqref{dft} can be rewritten by
letting $\ell  = \ell_1 n_2 + \ell_2$ and $k = k_1 + k_2 n_1$.  We then have:
\begin{equation}
\dfty[k_1 + k_2 n_1] = \label{eq:ct}\\
\sum_{\ell_2 = 0}^{n_2-1}
\left[\left(
    \sum_{\ell_1 = 0}^{n_1-1}
    \dftx[\ell_1 n_2 + \ell_2]
    \rootunity_{n_1}^{\ell_1 k_1}
  \right)
  \rootunity_{n}^{\ell_2 k_1}
\right]
\rootunity_{n_2}^{\ell_2 k_2}
\ , 
\end{equation}
where $k_{1,2} = 0,\ldots,n_{1,2} - 1$.  Thus, the algorithm computes
$n_2$ {\DFTs} of size $n_1$ (the inner sum), multiplies the result by
the so-called~\cite{GenSan66} \define{twiddle factors} $
\rootunity_{n}^{\ell_2 k_1} $, and finally computes $n_1$ {\DFTs} of size
$n_2$ (the outer sum).  This decomposition is then continued
recursively.  The literature uses the term \define{radix} to describe
an $n_1$~or~$n_2$ that is bounded (often constant); the small {\DFT}
of the radix is traditionally called a \define{butterfly}.

Many well-known variations are distinguished by the radix alone.  A
\define{decimation in time} (\define{{\DIT}}) algorithm uses $n_2$ as
the radix, while a \define{decimation in frequency} (\define{{\DIF}})
algorithm uses $n_1$ as the radix.  If multiple radices are used,
e.g. for $n$ composite but not a prime power, the algorithm is called
\define{mixed radix}.  A peculiar blending of radix~2 and 4 is called
\define{split radix}, which was proposed to minimize the count of
arithmetic
operations~\cite{Yavne68,Duhamel84,Vetterli84,Martens84,DuhamelVe90}
although it has been superseded in this
regard~\cite{Johnson07,Lundy07}.  {\FFTW} implements both {\DIT} and
{\DIF}, is mixed-radix with radices that are \emph{adapted} to the
hardware, and often uses much larger radices (e.g. radix~32) than were
once common.  On the other end of the scale, a ``radix'' of roughly
$\sqrt{n}$ has been called a \define{four-step} {\FFT} algorithm (or
\define{six-step}, depending on how many transposes one
performs)~\cite{Bailey90}; see \secref{cache} for some theoretical and
practical discussion of this algorithm.

A key difficulty in implementing the Cooley--Tukey {\FFT} is that the
$n_1$ dimension corresponds to discontiguous inputs $\ell_1$ in $\dftx$ but
contiguous outputs $k_1$ in $\dfty$, and vice-versa for $n_2$. This is a
matrix transpose for a single decomposition stage, and the composition
of all such transpositions is a (mixed-base) digit-reversal
permutation (or \define{bit-reversal}, for radix~2).  The resulting
necessity of discontiguous memory access and data re-ordering hinders
efficient use of hierarchical memory architectures (e.g., caches), so
that the optimal execution order of an {\FFT} for given hardware is
non-obvious, and various approaches have been proposed.

\begin{figure}[t]
\centering\includegraphics[width=0.70\columnwidth]{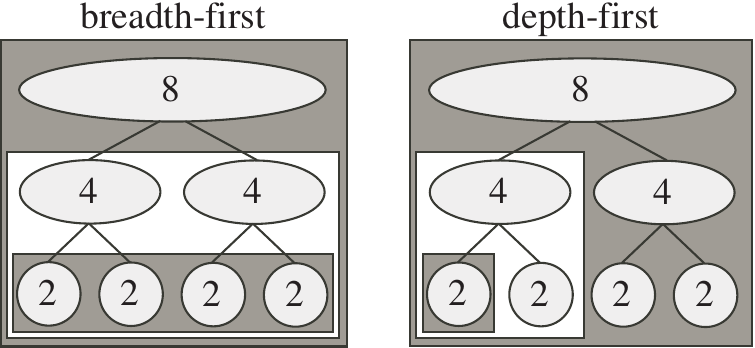}
\caption{Schematic of traditional breadth-first (left) vs. recursive
depth-first (right) ordering for radix-2 {\FFT} of size~8: the
computations for each nested box are completed before doing anything else in the
surrounding box.  Breadth-first computation performs all butterflies
of a given size at once, while depth-first computation completes one
subtransform entirely before moving on to the next (as in
\algref{radix2-rec}).}
\label{fig:breadthdepth}
\end{figure}

One ordering distinction is between recursion and iteration.  As
expressed above, the Cooley--Tukey algorithm could be thought of as
defining a tree of smaller and smaller {\DFTs}, as depicted in \figref{breadthdepth};
for example, a textbook radix-2 algorithm would divide size $n$ into
two transforms of size $n/2$, which are divided into four transforms
of size $n/4$, and so on until a base case is reached (in principle,
size 1).  This might naturally suggest a recursive implementation in
which the tree is traversed ``depth-first'' as in \figref{breadthdepth}(right) and
\algref{radix2-rec}---one size $n/2$ transform is solved completely
before processing the other one, and so on.  However, most traditional
{\FFT} implementations are non-recursive (with rare
exceptions~\cite{Singleton67}) and traverse the tree
``breadth-first''~\cite{VanLoan92} as in \figref{breadthdepth}(left)---in the radix-2
example, they would perform $n$ (trivial) size-1 transforms, then
$n/2$ combinations into size-2 transforms, then $n/4$ combinations
into size-4 transforms, and so on, thus making $\log_2 n$ passes over
the whole array.  In contrast, as we discuss in
\secref{struct:discussion}, {\FFTW} employs an explicitly
recursive strategy that encompasses \emph{both} depth-first and
breadth-first styles, favoring the former since it has some
theoretical and practical advantages as discussed in \secref{cache}.

\begin{algorithm}[t]
\begin{algorithmic}
\item[\textbf{function}] $\dfty[0,\ldots,n-1] \leftarrow \operatorname{\mathbf{recfft2}}(n, \dftx, \is)$:

\IF{$n = 1$}
\STATE $\dfty[0] \leftarrow \dftx[0]$
\ELSE
\STATE $\dfty[0,\ldots,n/2-1] \leftarrow \operatorname{\mathbf{recfft2}}(n/2, \dftx, 2\is)$
\STATE $\dfty[n/2, \ldots, n-1] \leftarrow \operatorname{\mathbf{recfft2}}(n/2, \dftx+\is, 2\is)$
\FOR{$k_1 = 0$ to $(n/2) - 1$}
\STATE $t \leftarrow \dfty[k_1]$
\STATE $\dfty[k_1] \leftarrow t + \rootunity_n^{k_1} \dfty[k_1 + n/2]$
\STATE $\dfty[k_1 + n/2] \leftarrow t - \rootunity_n^{k_1} \dfty[k_1 + n/2]$
\ENDFOR
\ENDIF

\end{algorithmic}
\caption{A depth-first recursive radix-2 {\DIT} Cooley--Tukey {\FFT} to
compute a {\DFT} of a power-of-two size $n=2^m$.  The input is an array
$\dftx$ of length $n$ with stride $\is$ (i.e., the inputs are $\dftx[\ell\is]$
for $\ell=0,\ldots,n-1$) and the output is an array $\dfty$ of length $n$ (with
stride 1), containing the {\DFT} of $\dftx$ [\eqref{dft}].  $\dftx+\is$
denotes the array beginning with $\dftx[\is]$.  This algorithm operates
out-of-place, produces in-order output, and does not require a
separate bit-reversal stage.
\label{alg:radix2-rec}}
\end{algorithm}

A second ordering distinction lies in how the digit-reversal is
performed.  The classic approach is a single, separate digit-reversal
pass following or preceding the arithmetic computations; this approach
is so common and so deeply embedded into {\FFT} lore that many
practitioners find it difficult to imagine an {\FFT} without an
explicit bit-reversal stage.  Although this pass requires only $O(n)$
time~\cite{Karp96}, it can still be non-negligible, especially if the
data is out-of-cache; moreover, it neglects the possibility that data
reordering during the transform may improve memory locality.  Perhaps
the oldest alternative is the Stockham \define{auto-sort}
{\FFT}~\cite{Stockham66,VanLoan92}, which transforms back and forth
between two arrays with each butterfly, transposing one digit each
time, and was popular to improve contiguity of access for vector
computers~\cite{Swarztrauber82}.  Alternatively, an explicitly
recursive style, as in {\FFTW}, performs the digit-reversal implicitly
at the ``leaves'' of its computation when operating out-of-place
(\secref{struct:discussion}).  A simple example of this
style, which computes in-order output using an out-of-place radix-2
{\FFT} without explicit bit-reversal, is shown in \algref{radix2-rec}
[corresponding to \figref{breadthdepth}(right)].  To operate in-place
with $O(1)$ scratch storage, one can interleave small matrix
transpositions with the
butterflies~\cite{JohnsonBu84,Temperton91,Qian94,Hegland94}, and a
related strategy in {\FFTW}~\cite{FFTW05} is briefly described by
\secref{struct:discussion}.

Finally, we should mention that there are many {\FFTs} entirely
distinct from Cooley--Tukey.  Three notable such algorithms are the
\define{prime-factor algorithm} for $\gcd(n_1,n_2)=1$~\cite[page
619]{OppenheimSha89}, along with Rader's~\cite{Rader68} and
Bluestein's~\cite{Bluestein68,RabinerScRa69,OppenheimSha89} algorithms
for prime~$n$.  {\FFTW} implements the first two in its codelet
generator for hard-coded $n$ (\secref{genfft}) and the latter two for
general prime~$n$ (\secreftwo{struct:prime}{generality}).  There is also the Winograd
{\FFT}~\cite{Winograd78,Heideman86,Duhamel90,DuhamelVe90}, which
minimizes the number of multiplications at the expense of a large
number of additions; this trade-off is not beneficial on current
processors that have specialized hardware multipliers.

\section{{\FFTs} and the Memory Hierarchy}
\label{sec:cache}

There are many complexities of computer architectures that impact the
optimization of {\FFT} implementations, but one of the most pervasive
is the memory hierarchy.  On any modern general-purpose computer,
memory is arranged into a hierarchy of storage devices with increasing
size and decreasing speed: the fastest and smallest memory being the
CPU registers, then two or three levels of cache, then the main-memory
RAM, then external storage such as hard disks.\footnote{A hard disk is
utilized by ``out-of-core'' {\FFT} algorithms for very large $n$~\cite{VanLoan92}, but these algorithms appear to have been largely superseded in practice
by both the gigabytes of memory now common on personal computers and,
for extremely large $n$, by algorithms for distributed-memory parallel
computers.}  Most of these levels are managed automatically by
the hardware to hold the most-recently-used data from the next level
in the hierarchy.\footnote{This includes the registers: on current
``x86'' processors, the user-visible instruction set (with a small
number of floating-point registers) is internally translated at
runtime to RISC-like ``$\mu$-ops'' with a much larger number of
physical \define{rename registers} that are allocated automatically.}
There are many complications, however, such as limited cache
associativity (which means that certain locations in memory cannot be
cached simultaneously) and cache lines (which optimize the cache for
contiguous memory access), which are reviewed in numerous textbooks on
computer architectures [refs].  In this section, we focus on the
simplest abstract principles of memory hierarchies in order to grasp
their fundamental impact on {\FFTs}.

Because access to memory is in many cases the slowest part of the
computer, especially compared to arithmetic, one wishes to load as
much data as possible in to the faster levels of the hierarchy, and
then perform as much computation as possible before going back to the
slower memory devices.  This is called \define{temporal locality}
[ref]: if a given datum is used more than once, we arrange the
computation so that these usages occur as close together as possible
in time.

\subsection{Understanding {\FFTs} with an ideal cache}
\label{sec:cache:ideal}

To understand temporal-locality strategies at a basic level, in this
section we will employ an idealized model of a cache in a two-level
memory hierarchy, as defined in \citeasnoun{FrigoLe99}.  This
\define{ideal cache} stores $Z$ data items from main memory
(e.g. complex numbers for our purposes): when the processor loads a
datum from memory, the access is quick if the datum is already in the
cache (a \define{cache hit}) and slow otherwise (a \define{cache
miss}, which requires the datum to be fetched into the cache).  When a
datum is loaded into the cache,\footnote{More generally, one can
assume that a \define{cache line} of $L$ consecutive data items are
loaded into the cache at once, in order to exploit spatial locality.
The ideal-cache model in this case requires that the cache be
\define{tall}: $Z = \Omega(L^2)$~\cite{FrigoLe99}.} it must replace
some other datum, and the ideal-cache model assumes that the optimal
replacement strategy is used~\cite{Belady66}: the new datum replaces
the datum that will not be needed for the longest time in the future;
in practice, this can be simulated to within a factor of two by
replacing the least-recently used datum~\cite{FrigoLe99}, but ideal
replacement is much simpler to analyze.  Armed with this ideal-cache
model, we can now understand some basic features of {\FFT}
implementations that remain essentially true even on real cache
architectures.  In particular, we want to know the \define{cache
complexity}, the number $Q(n;Z)$ of cache misses for an {\FFT} of size
$n$ with an ideal cache of size $Z$, and what algorithm choices reduce
this complexity.

First, consider a textbook radix-2 algorithm, which divides $n$ by~$2$
at each stage and operates breadth-first as in
\figref{breadthdepth}(left), performing all butterflies of a given
size at a time.  If $n > Z$, then each pass over the array incurs
$\Theta(n)$ cache misses to reload the data, and there are $\log_2 n$
passes, for $\Theta(n \log_2 n)$ cache misses in total---no temporal
locality at all is exploited!

One traditional solution to this problem is \define{blocking}: the
computation is divided into maximal blocks that fit into the cache,
and the computations for each block are completed before moving on to
the next block.  Here, a block of $Z$ numbers can fit into the
cache\footnote{Of course, $O(n)$ additional storage may be required
for twiddle factors, the output data (if the {\FFT} is not in-place),
and so on, but these only affect the $n$ that fits into cache by a
constant factor and hence do not impact cache-complexity analysis.  We
won't worry about such constant factors in this section.}  (not
including storage for twiddle factors and so on), and thus the natural
unit of computation is a sub-{\FFT} of size $Z$.  Since each of these
blocks involves $\Theta(Z \log Z)$ arithmetic operations, and there
are $\Theta(n \log n)$ operations overall, there must be
$\Theta(\frac{n}{Z} \log_Z n)$ such blocks.  More explicitly, one
could use a radix-$Z$ Cooley--Tukey algorithm, breaking $n$ down by
factors of $Z$ [or $\Theta(Z)$] until a size $Z$ is reached: each
stage requires $n/Z$ blocks, and there are $\log_Z n$ stages, again
giving $\Theta(\frac{n}{Z} \log_Z n)$ blocks overall.  Since each
block requires $Z$ cache misses to load it into cache, the cache
complexity $Q_b$ of such a blocked algorithm is
\begin{equation}
\label{eq:cache-opt}
Q_b(n;Z) = \Theta(n \log_Z n) .
\end{equation}
In fact, this complexity is rigorously \emph{optimal} for Cooley--Tukey
{\FFT} algorithms~\cite{HongKu81}, and immediately points us towards
\emph{large radices} (not radix~2!) to exploit caches effectively in
{\FFTs}.

However, there is one shortcoming of any blocked {\FFT} algorithm: it
is \define{cache aware}, meaning that the implementation depends
explicitly on the cache size $Z$.  The implementation must be modified
(e.g. changing the radix) to adapt to different machines as the cache
size changes.  Worse, as mentioned above, actual machines have
multiple levels of cache, and to exploit these one must perform
multiple levels of blocking, each parameterized by the corresponding
cache size.  In the above example, if there were a smaller and faster
cache of size $z < Z$, the size-$Z$ sub-{\FFTs} should themselves be
performed via radix-$z$ Cooley--Tukey using blocks of size $z$.  And so
on.  There are two paths out of these difficulties: one is
self-optimization, where the implementation automatically adapts
itself to the hardware (implicitly including any cache sizes), as
described in \secref{struct}; the other is to exploit
\define{cache-oblivious} algorithms.  {\FFTW} employs both of these
techniques.

The goal of cache-obliviousness is to structure the algorithm so that
it exploits the cache without having the cache size as a parameter:
the same code achieves the same asymptotic cache complexity regardless
of the cache size $Z$.  An \define{optimal cache-oblivious} algorithm
achieves the \emph{optimal} cache complexity (that is, in an
asymptotic sense, ignoring constant factors).  Remarkably, optimal
cache-oblivious algorithms exist for many problems, such as matrix
multiplication, sorting, transposition, and {\FFTs}~\cite{FrigoLe99}.
Not all cache-oblivious algorithms are optimal, of course---for
example, the textbook radix-2 algorithm discussed above is
``pessimal'' cache-oblivious (its cache complexity is independent of
$Z$ because it always achieves the worst case!).

For instance, \figref{breadthdepth}(right) and \algref{radix2-rec} shows a way to
obliviously exploit the cache with a radix-2 Cooley--Tukey algorithm,
by ordering the computation depth-first rather than breadth-first.
That is, the {\DFT} of size $n$ is divided into two {\DFTs} of size
$n/2$, and one {\DFT} of size $n/2$ is \emph{completely finished}
before doing \emph{any} computations for the second {\DFT} of size
$n/2$.  The two subtransforms are then combined using $n/2$ radix-2
butterflies, which requires a pass over the array and (hence $n$ cache
misses if $n>Z$). This process is repeated recursively until a
base-case (e.g. size $2$) is reached.  The cache complexity $Q_2(n;Z)$
of this algorithm satisfies the recurrence
\begin{equation}
\label{eq:cache-radix2-rec}
Q_2(n;Z) = \begin{cases} n & n \leq Z \\ 
                     2Q_2(n/2;Z) + \Theta(n) & \mathrm{otherwise} \end{cases} .
\end{equation}
The key property is this: once the recursion reaches a size $n \leq
Z$, the subtransform fits into the cache and no further misses are
incurred.  The algorithm does not ``know'' this and continues
subdividing the problem, of course, but all of those further
subdivisions are in-cache because they are performed in the same
depth-first branch of the tree.  The solution of
\eqref{cache-radix2-rec} is
\begin{equation}
\label{eq:cache-subopt}
Q_2(n;Z) = \Theta(n \log [n/Z]) .
\end{equation}
This is worse than the theoretical optimum $Q_b(n;Z)$ from
\eqref{cache-opt}, but it is cache-oblivious ($Z$ never entered the
algorithm) and exploits at least \emph{some} temporal
locality.\footnote{This advantage of depth-first recursive
implementation of the radix-2 {\FFT} was pointed out many years ago
by Singleton (where the ``cache'' was core
memory)~\cite{Singleton67}.}  On the other hand, when it is combined
with {\FFTW}'s self-optimization and larger radices in
\secref{struct}, this algorithm actually performs very well until $n$
becomes extremely large.  By itself, however, \algref{radix2-rec} must
be modified to attain adequate performance for reasons that have
nothing to do with the cache.  These practical issues are discussed
further in \secref{cache:practice}.

There exists a different recursive {\FFT} that is \emph{optimal}
cache-oblivious, however, and that is the radix-$\sqrt{n}$
``four-step'' Cooley--Tukey algorithm (again executed recursively,
depth-first)~\cite{FrigoLe99}.  The cache complexity $Q_o$ of this algorithm
satisfies the recurrence:
\begin{equation}
\label{eq:cache-4step}
Q_o(n;Z) = \begin{cases} n & n \leq Z \\ 
                     2 \sqrt{n} Q_o(\sqrt{n};Z) + \Theta(n) & \mathrm{otherwise}
         \end{cases} .
\end{equation}
That is, at each stage one performs $\sqrt{n}$ {\DFTs} of size
$\sqrt{n}$ (recursively), then multiplies by the $\Theta(n)$ twiddle
factors (and does a matrix transposition to obtain in-order output),
then finally performs another $\sqrt{n}$ {\DFTs} of size $\sqrt{n}$.  The
solution of \eqref{cache-4step} is $Q_o(n;Z) = \Theta(n \log_Z n)$, the
same as the optimal cache complexity \eqref{cache-opt}!

These algorithms illustrate the basic features of most optimal
cache-oblivious algorithms: they employ a recursive divide-and-conquer
strategy to subdivide the problem until it fits into cache, at which
point the subdivision continues but no further cache misses are
required.  Moreover, a cache-oblivious algorithm exploits all levels
of the cache in the same way, so an optimal cache-oblivious algorithm
exploits a multi-level cache optimally as well as a two-level
cache~\cite{FrigoLe99}: the multi-level ``blocking'' is implicit in
the recursion.

\subsection{Cache-obliviousness in practice}
\label{sec:cache:practice}

Even though the radix-$\sqrt{n}$ algorithm is optimal cache-oblivious,
it does not follow that {\FFT} implementation is a solved problem. The
optimality is only in an asymptotic sense, ignoring constant factors,
$O(n)$ terms, etcetera, all of which can matter a great deal in
practice.  For small or moderate $n$, quite different algorithms may
be superior, as discussed in \secref{cache:fftw}.  Moreover, real
caches are inferior to an ideal cache in several ways.  The
unsurprising consequence of all this is that cache-obliviousness, like
any complexity-based algorithm property, does not absolve one from the
ordinary process of software optimization.  At best, it reduces the
amount of memory/cache tuning that one needs to perform, structuring
the implementation to make further optimization easier and more
portable.

Perhaps most importantly, one needs to perform an optimization that
has almost nothing to do with the caches: the recursion must be
``coarsened'' to amortize the function-call overhead and to enable
compiler optimization.  For example, the simple pedagogical code of
\algref{radix2-rec} recurses all the way down to $n=1$, and hence
there are $\approx 2n$ function calls in total, so that every data
point incurs a two-function-call overhead on average.  Moreover, the
compiler cannot fully exploit the large register sets and instruction-level
parallelism of modern processors with an $n=1$ function
body.\footnote{In principle, it might be possible for a compiler to
automatically coarsen the recursion, similar to how compilers can
partially unroll loops.  We are currently unaware of any
general-purpose compiler that performs this optimization, however.}
These problems can be effectively erased, however, simply by making
the base cases larger, e.g. the recursion could stop when $n=32$ is
reached, at which point a highly optimized hard-coded {\FFT} of that
size would be executed.  In {\FFTW}, we produced this sort of large
base-case using a specialized code-generation program described in
\secref{genfft}.

One might get the impression that there is a strict dichotomy that
divides cache-aware and cache-oblivious algorithms, but the two are
not mutually exclusive in practice.  Given an implementation of a
cache-oblivious strategy, one can further optimize it for the cache
characteristics of a particular machine in order to improve the
constant factors.  For example, one can tune the radices used, the
transition point between the radix-$\sqrt{n}$ algorithm and the
bounded-radix algorithm, or other algorithmic choices as described
in \secref{cache:fftw}.  The advantage of starting cache-aware tuning with a
cache-oblivious approach is that the starting point already exploits
all levels of the cache to some extent, and one has reason to hope
that good performance on one machine will be more portable to other
architectures than for a purely cache-aware ``blocking'' approach.  In
practice, we have found this combination to be very successful with
{\FFTW}.

\subsection{Memory strategies in {\FFTW}}
\label{sec:cache:fftw}

The recursive cache-oblivious strategies described above form a useful
starting point, but {\FFTW} supplements them with a number of
additional tricks, and also exploits cache-obliviousness in
less-obvious forms.

We currently find that the general radix-$\sqrt{n}$ algorithm is
beneficial only when $n$ becomes very large, on the order of $2^{20}
\approx 10^6$.  In practice, this means that we use at most a single
step of radix-$\sqrt{n}$ (two steps would only be used for $n \gtrsim
2^{40}$).  The reason for this is that the implementation of
radix~$\sqrt{n}$ is less efficient than for a bounded radix: the
latter has the advantage that an entire radix butterfly can be
performed in hard-coded loop-free code within local
variables/registers, including the necessary permutations and twiddle
factors.  

Thus, for more moderate $n$, {\FFTW} uses depth-first recursion with a
bounded radix, similar in spirit to \algref{radix2-rec} but with much
larger radices (radix~32 is common) and base cases (size 32 or 64 is
common) as produced by the code generator of \secref{genfft}.  The
self-optimization described in \secref{struct} allows the choice of
radix and the transition to the radix-$\sqrt{n}$ algorithm to be tuned
in a cache-aware (but entirely automatic) fashion.

For small $n$ (including the radix butterflies and the base cases of
the recursion), hard-coded {\FFTs} ({\FFTW}'s \define{codelets}) are
employed.  However, this gives rise to an interesting problem: a
codelet for (e.g.) $n=64$ is $\sim 2000$ lines long, with hundreds of
variables and over 1000 arithmetic operations that can be executed in
many orders, so what order should be chosen?  The key problem here is
the efficient use of the CPU registers, which essentially form a
nearly ideal, fully associative cache.  Normally, one relies on the
compiler for all code scheduling and register allocation, but but the
compiler needs help with such long blocks of code (indeed, the general
register-allocation problem is NP-complete).  In particular, {\FFTW}'s
generator knows more about the code than the compiler---the generator
knows it is an FFT, and therefore it can use an optimal
cache-oblivious schedule (analogous to the radix-$\sqrt{n}$ algorithm)
to order the code independent of the number of
registers~\cite{Frigo99}.  The compiler is then used only for local
``cache-aware'' tuning (both for register allocation and the CPU
pipeline).\footnote{One practical difficulty is that some
``optimizing'' compilers will tend to greatly re-order the code,
destroying {\FFTW}'s optimal schedule.  With GNU gcc, we circumvent this
problem by using compiler flags that explicitly disable certain stages of the
optimizer.}  As a practical matter, one consequence of this scheduler
is that {\FFTW}'s machine-independent codelets are no slower than
machine-specific codelets generated by an automated search and
optimization over many possible codelet implementations, as performed
by the {\SPIRAL} project \cite[Figure~3]{XiongPa01}.

(When implementing hard-coded base cases, there is another choice because a
loop of small transforms is always required.  Is it better to
implement a hard-coded {\FFT} of size 64, for example, or an unrolled
loop of four size-16 {\FFTs}, both of which operate on the same amount
of data?  The former should be more efficient because it performs more
computations with the same amount of data, thanks to the $\log n$
factor in the {\FFT}'s $n \log n$ complexity.)

In addition, there are many other techniques that {\FFTW} employs to
supplement the basic recursive strategy, mainly to address the fact
that cache implementations strongly favor accessing consecutive
data---thanks to cache lines, limited associativity, and direct
mapping using low-order address bits (accessing data at power-of-two
intervals in memory, which is distressingly common in {\FFTs}, is thus
especially prone to cache-line conflicts).  Unfortunately, the known
{\FFT} algorithms inherently involve some non-consecutive access
(whether mixed with the computation or in separate
bit-reversal/transposition stages).  There are many optimizations in
{\FFTW} to address this.  For example, the data for several
butterflies at a time can be copied to a small buffer before computing
and then copied back, where the copies and computations involve more
consecutive access than doing the computation directly in-place.  Or,
the input data for the subtransform can be copied from
(discontiguous) input to (contiguous) output before performing the
subtransform in-place (see \secref{struct:indirect}), rather than
performing the subtransform directly out-of-place (as in
\algref{radix2-rec}). Or, the order of loops can be interchanged in
order to push the outermost loop from the first radix step [the $\ell_2$
loop in \eqref{ct}] down to the leaves, in order to make the input
access more consecutive (see \secref{struct:discussion}).  Or, the twiddle
factors can be computed using a smaller look-up table (fewer memory
loads) at the cost of more arithmetic (see \secref{accuracy}).  The
choice of whether to use any of these techniques, which come into play
mainly for moderate $n$ ($2^{13} < n < 2^{20}$), is made by the
self-optimizing planner as described in the next section.

\section{Adaptive Composition of {\FFT} Algorithms}
\label{sec:struct}

As alluded to several times already, {\FFTW} implements a wide variety
of {\FFT} algorithms (mostly rearrangements of Cooley--Tukey) and
selects the ``best'' algorithm for a given $n$ automatically.  In this
section, we describe how such self-optimization is implemented, and
especially how {\FFTW}'s algorithms are structured as a composition of
algorithmic fragments.  These techniques in {\FFTW} are described in
greater detail elsewhere~\cite{FFTW05}, so here we will focus only on
the essential ideas and the motivations behind them.

An {\FFT} algorithm in {\FFTW} is a composition of algorithmic steps
called a \define{plan}.  The algorithmic steps each solve a certain
class of \define{problems} (either solving the problem directly or
recursively breaking it into sub-problems of the same type).  The
choice of plan for a given problem is determined by a \define{planner}
that selects a composition of steps, either by runtime measurements to
pick the fastest algorithm, or by heuristics, or by loading a
pre-computed plan.  These three pieces: problems, algorithmic steps,
and the planner, are discussed in the following subsections.

\subsection{The problem to be solved}

In early versions of {\FFTW}, the only choice made by the planner was
the sequence of radices~\cite{FrigoJo98}, and so each step of the plan
took a {\DFT} of a given size $n$, possibly with discontiguous
input/output, and reduced it (via a radix $r$) to {\DFTs} of size
$n/r$, which were solved recursively.  That is, each step solved the
following problem: given a size $n$, an \define{input pointer}
$\dfti$, an \define{input stride} $\is$, an \define{output pointer}
$\dfto$, and an \define{output stride} $\os$, it computed the {\DFT}
of $\dfti[\ell\is]$ for $0 \leq \ell < n$ and stored the result in
$\dfto[k\os]$ for $0 \leq k < n$.  However, we soon found that we
could not easily express many interesting algorithms within this
framework; for example, \define{in-place} ($\dfti = \dfto$) {\FFTs}
that do not require a separate bit-reversal
stage~\cite{JohnsonBu84,Temperton91,Qian94,Hegland94}.  It became
clear that the key issue was not the choice of algorithms, as we had
first supposed, but the definition of the problem to be solved.
Because only problems that can be expressed can be solved, the
representation of a problem determines an outer bound to the space of
plans that the planner can explore, and therefore it ultimately
constrains {\FFTW}'s performance.

The difficulty with our initial $(n,\dfti,\is,\dfto,\os)$ problem
definition was that it forced each algorithmic step to address only a
single {\DFT}.  In fact, {\FFTs} break down {\DFTs} into
\emph{multiple} smaller {\DFTs}, and it is the \emph{combination} of
these smaller transforms that is best addressed by many algorithmic
choices, especially to rearrange the order of memory accesses between
the subtransforms.  Therefore, we redefined our notion of a problem
in {\FFTW} to be not a single {\DFT}, but rather a \emph{loop} of
{\DFTs}, and in fact \emph{multiple nested loops} of {\DFTs}.  The
following sections describe some of the new algorithmic steps that
such a problem definition enables, but first we will define the problem
more precisely.

{\DFT} problems in {\FFTW} are expressed in terms of structures called
{\IOTENSs},\footnote{{\IOTENSs} are unrelated to the tensor-product
notation used by some other authors to describe {\FFT}
algorithms~\cite{VanLoan92,Puschel05}.}  which in turn are described
in terms of ancillary structures called {\IODIMs}.  An \define{\IODIM}
$\aiodim$ is a triple $\aiodim = \iodim{\n}{\is}{\os}$, where $\n$~is
a non-negative integer called the \define{length}, $\is$ is an integer
called the \define{input stride}, and $\os$ is an integer called the
\define{output stride}.  An \define{\IOTENS}
$\aiotens=\iotens{\aiodim_1, \aiodim_2, \ldots, \aiodim_\rank}$ is a
set of {\IODIMs}.  The non-negative integer~$\rank = \rankof{\aiotens}$
is called the \define{rank} of the {\IOTENS}.  A \define{{\DFT}
problem}, denoted by $\dftprob{\dftn}{\dftv}{\dfti}{\dfto}$, consists
of two {\IOTENSs} $\dftn$ and $\dftv$, and of two \define{pointers}
$\dfti$ and $\dfto$.  Informally, this describes $\rankof{\dftv}$
nested loops of $\rankof{\dftn}$-dimensional {\DFTs} with input data
starting at memory location~$\dfti$ and output data starting
at~$\dfto$.

For simplicity, let us consider only one-dimensional {\DFTs}, so that
$\dftn = \iotens{\iodim{\n}{\is}{\os}}$ implies a $\DFT$ of length $n$
on input data with stride $\is$ and output data with stride $\os$,
much like in the original {\FFTW} as described above.  The main new
feature is then the addition of zero or more ``loops'' $\dftv$.  More
formally, $\dftprob{\dftn}{\iotens{\iodim{\n}{\is}{\os}} \union
\dftv}{\dfti}{\dfto}$ is recursively defined as a ``loop'' of $\n$
problems: for all $0 \leq k < \n$, do all computations in
$\dftprob{\dftn}{\dftv}{\dfti + k \cdot \is}{\dfto + k \cdot \os}$.
The case of multi-dimensional {\DFTs} is defined more precisely
elsewhere~\cite{FFTW05}, but essentially each {\IODIM} in $\dftn$ gives one
dimension of the transform.

We call $\dftn$ the \define{size} of the problem.  The \define{rank}
of a problem is defined to be the rank of its size (i.e., the
dimensionality of the {\DFT}).  Similarly, we call $\dftv$ the
\define{vector size} of the problem, and the \define{vector rank} of a
problem is correspondingly defined to be the rank of its vector size.
Intuitively, the vector size can be interpreted as a set of ``loops''
wrapped around a single {\DFT}, and we therefore refer to a single
{\IODIM} of $\dftv$ as a \define{vector loop}.  (Alternatively, one
can view the problem as describing a {\DFT} over a
$\rankof{\dftv}$-dimensional vector space.)  The problem does not
specify the order of execution of these loops, however, and therefore
{\FFTW} is free to choose the fastest or most convenient order.

\subsubsection{{\DFT} problem examples}

A more detailed discussion of the space of problems in {\FFTW} can be
found in \citeasnoun{FFTW05}, but a simple understanding can be gained
by examining a few examples demonstrating that the {\IOTENS}
representation is sufficiently general to cover many situations that
arise in practice, including some that are not usually considered to
be instances of the {\DFT}.

A single one-dimensional {\DFT} of length $n$, with
stride-1 input $\dftx$ and output $\dfty$, as in \eqref{dft}, is
denoted by the problem
$\dftprob{\iotens{\iodim{n}{1}{1}}}{\iotensone}{\dftx}{\dfty}$ (no
loops: vector-rank zero).

As a more complicated example, suppose we have an $n_1 \times n_2$
matrix $\dftx$ stored as $n_1$ consecutive blocks of contiguous
length-$n_2$ rows (this is called \emph{row-major} format).  The
in-place {\DFT} of all the \emph{rows} of this matrix would be denoted
by the problem
$\dftprob{\iotens{\iodim{n_2}{1}{1}}}{\iotens{\iodim{n_1}{n_2}{n_2}}}{\dftx}{\dftx}$:
a length-$n_1$ loop of size-$n_2$ contiguous {\DFTs}, where each
iteration of the loop offsets its input/output data by a stride $n_2$.
Conversely, the in-place {\DFT} of all the \emph{columns} of this
matrix would be denoted by
$\dftprob{\iotens{\iodim{n_1}{n_2}{n_2}}}{\iotens{\iodim{n_2}{1}{1}}}{\dftx}{\dftx}$:
compared to the previous example, $\dftn$ and $\dftv$ are swapped.  In
the latter case, each {\DFT} operates on discontiguous data, and
{\FFTW} might well choose to interchange the loops: instead of
performing a loop of {\DFTs} computed individually, the subtransforms
themselves could act on $n_1$-component vectors, as described in 
\secref{struct:plans}.

A size-1 {\DFT} is simply a copy $Y[0] = X[0]$, and here this can also
be denoted by $\dftn = \iotensone$ (rank zero, a ``zero-dimensional'' {\DFT}).
This allows {\FFTW}'s problems to represent many kinds of copies and
permutations of the data within the same problem framework, which is
convenient because these sorts of operations arise frequently in
{\FFT} algorithms.  For example, to copy $n$ consecutive numbers from
$\dfti$ to $\dfto$, one would use the rank-zero problem
$\dftprob{\iotensone}{\iotens{\iodim{n}{1}{1}}}{\dfti}{\dfto}$.  More
interestingly, the in-place \emph{transpose} of an $n_1 \times n_2$
matrix $\dftx$ stored in row-major format, as described above, is
denoted by $\dftprob{\iotensone}{\iotens{\iodim{n_1}{n_2}{1},
\iodim{n_2}{1}{n_1}}}{\dftx}{\dftx}$ (rank zero, vector-rank two).

\subsection{The space of plans in {\FFTW}}
\label{sec:struct:plans}

Here, we describe a subset of the possible plans considered by
{\FFTW}; while not exhaustive~\cite{FFTW05}, this subset is enough to
illustrate the basic structure of {\FFTW} and the necessity of
including the vector loop(s) in the problem definition to enable
several interesting algorithms.  The plans that we now describe
usually perform some simple ``atomic'' operation, and it may not be
apparent how these operations fit together to actually compute
{\DFTs}, or why certain operations are useful at all.  We shall
discuss those matters in~\secref{struct:discussion}.

Roughly speaking, to solve a general {\DFT} problem, one must perform
three tasks.  First, one must reduce a problem of arbitrary vector
rank to a set of loops nested around a problem of vector rank~0, i.e.,
a single (possibly multi-dimensional) {\DFT}.  Second, one must reduce
the multi-dimensional {\DFT} to a sequence of of rank-1 problems,
i.e., one-dimensional {\DFTs}; for simplicity, however, we do not
consider multi-dimensional {\DFTs} below.  Third, one must solve the
rank-1, vector rank-0 problem by means of some {\DFT} algorithm such
as Cooley--Tukey.  These three steps need not be executed in the stated
order, however, and in fact, almost every permutation and interleaving
of these three steps leads to a correct {\DFT} plan. The choice of the
set of plans explored by the planner is critical for the usability of
the {\FFTW} system: the set must be large enough to contain the
fastest possible plans, but it must be small enough to keep the
planning time acceptable.

\subsubsection{Rank-0 plans}
\label{sec:struct:rank0} 
The rank-0 problem $\dftprob{\iotensone}{\dftv}{\dfti}{\dfto}$ denotes
a permutation of the input array into the output array.  {\FFTWthree}
does not solve arbitrary rank-0 problems, only the following two
special cases that arise in practice.

\begin{itemize}
\item When $\rankof{\dftv} = 1$ and $\dfti \neq \dfto$, {\FFTWthree}
  produces a plan that copies the input array into the output array.
  Depending on the strides, the plan consists of a
  loop or, possibly, of a call to the ANSI~C function \texttt{memcpy}, which is
  specialized to copy contiguous regions of memory.
\item When $\rankof{\dftv} = 2$, $\dfti = \dfto$, and the strides
  denote a matrix-transposition problem, {\FFTWthree} creates a plan
  that transposes the array in-place.  {\FFTWthree} implements the
  square transposition
  $\dftprob{\iotensone}{\iotens{\iodim{\n}{\is}{\os},
  \iodim{\n}{\os}{\is}}} {\dfti}{\dfto}$ by means of the
  ``cache-oblivious'' algorithm from \cite{FrigoLe99}, which is fast
  and, in theory, uses the cache optimally regardless of the cache
  size.  A generalization of this idea is employed for non-square
  transpositions with a large common factor or a small difference
  between the dimensions, adapting algorithms from \citeasnoun{Dow95}.
\end{itemize}

\subsubsection{Rank-1 plans}
\label{sec:struct:rank1} 
Rank-1 {\DFT} problems denote ordinary one-dimensional Fourier
transforms.  {\FFTWthree} deals with most rank-1 problems as follows.

\paragraph{Direct plans}
\label{sec:struct:rank1:direct} 

When the {\DFT} rank-1 problem is ``small enough'' (usually, $n \leq
64$), {\FFTWthree} produces a \define{direct plan} that solves the
problem directly.  These plans operate by calling a fragment of C code
(a \define{codelet}) specialized to solve problems of one particular
size, whose generation is described in \secref{genfft}.  More
precisely, the codelets compute a loop ($\rankof{\dftv} \leq 1$) of
small {\DFTs}.

\paragraph{Cooley--Tukey plans}
\label{sec:struct:rank1:Cooley--Tukey} 

For problems of the form
$\dftprob{\iotens{\iodim{\n}{\is}{\os}}}{\dftv}{\dfti}{\dfto}$ where
$\n = \radix\m$, {\FFTWthree} generates a plan that implements a
radix-$\radix$ Cooley--Tukey algorithm (\secref{fftintro}).
Both decimation-in-time and decimation-in-frequency plans are supported, with both small fixed radices (usually, $\radix \leq 64$) produced by the codelet generator (\secref{genfft}) and also arbitrary radices (e.g. radix-$\sqrt{n}$).

The most common case is a \define{decimation in time} (\define{{\DIT}}) plan, corresponding to a
\define{radix} $\radix = n_2$ (and thus $\m = n_1$): it first solves
$\dftprob{\iotens{\iodim{\m}{\radix\cdot\is}{\os}}}
{\dftv\union\iotens{\iodim{\radix}{\is}{\m\cdot\os}}} {\dfti}{\dfto}$,
then multiplies the output array $\dfto$ by the twiddle factors, and
finally solves
$\dftprob{\iotens{\iodim{\radix}{\m\cdot\os}{\m\cdot\os}}}
{\dftv\union\iotens{\iodim{\m}{\os}{\os}}} {\dfto}{\dfto}$.  For
performance, the last two steps are not planned independently, but are
fused together in a single ``twiddle'' codelet---a fragment of C code
that multiplies its input by the twiddle factors and performs a {\DFT}
of size~$\radix$, operating in-place on $\dfto$.

\subsubsection{Plans for higher vector ranks}
\label{sec:struct:vrankr} 
These plans extract a vector loop to reduce a {\DFT} problem to a
problem of lower vector rank, which is then solved recursively.  Any
of the vector loops of $\dftv$ could be extracted in this way, leading
to a number of possible plans corresponding to different loop
orderings.

Formally, to solve $\dftprob{\dftn}{\dftv}{\dfti}{\dfto}$, where
$\dftv=\iotens{\iodim{\n}{\is}{\os}} \union \dftv_1$, {\FFTWthree}
generates a loop that, for all $\ndx$ such that $0\leq \ndx < \n$,
invokes a plan for
$\dftprob{\dftn}{\dftv_1}{\dfti+\ndx\cdot\is}{\dfto+\ndx\cdot\os}$.

\subsubsection{Indirect plans} 
\label{sec:struct:indirect}
Indirect plans transform a {\DFT} problem that requires some data
shuffling (or discontiguous operation) into a problem that requires no
shuffling plus a rank-0 problem that performs the shuffling.

Formally, to solve $\dftprob{\dftn}{\dftv}{\dfti}{\dfto}$ where
$\rankof{\dftn}>0$, {\FFTWthree} generates a plan that first solves
$\dftprob{\iotensone}{\dftn\union\dftv}{\dfti}{\dfto}$, and then
solves $\dftprob{\copyo{\dftn}}{\copyo{\dftv}}{\dfto}{\dfto}$.  Here
we define $\copyo{\aiotens}$ to be the {\IOTENS}
$\iotens{\iodim{\n}{\os}{\os} \suchthat
\iodim{\n}{\is}{\os}\in\aiotens}$: that is, it replaces the input
strides with the output strides.  Thus, an indirect plan first
rearranges/copies the data to the output, then solves the problem in
place.

\subsubsection{Plans for prime sizes}
\label{sec:struct:prime}
As discussed in \secref{generality}, it turns out to be surprisingly
useful to be able to handle large prime $n$ (or large prime factors).
\define{Rader plans} implement the algorithm from \citeasnoun{Rader68}
to compute one-dimensional {\DFTs} of prime size in $\Theta(n \log n)$
time.  \define{Bluestein plans} implement Bluestein's ``chirp-z''
algorithm, which can also handle prime $n$ in $\Theta(n \log n)$
time~\cite{Bluestein68,RabinerScRa69,OppenheimSha89}.  \define{Generic plans}
implement a naive $\Theta(n^2)$ algorithm (useful for $n \lesssim
100$).

\subsubsection{Discussion}
\label{sec:struct:discussion}

Although it may not be immediately apparent, the combination of the
recursive rules in \secref{struct:plans} can produce a number of
useful algorithms.  To illustrate these compositions, we discuss in
particular three issues: depth- vs. breadth-first, loop reordering,
and in-place transforms.

\begin{figure}
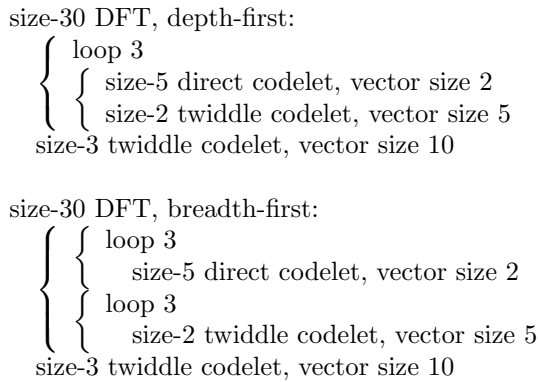

\[
\begin{array}{l}
 \begin{array}{l}
 \mbox{size-30 {\DFT}, depth-first:}\\
   \quad\left \{ 
     \begin{array}{l}
       \mbox{loop 3} \\
       \left \{
         \begin{array}{l}
           \mbox{size-5 direct codelet, vector size 2} \\
           \mbox{size-2 twiddle codelet, vector size 5} \\
         \end{array}
       \right . \\
    \end{array}
   \right .\\
   \quad\mbox{size-3 twiddle codelet, vector size 10}
 \end{array} \\
\\
 \begin{array}{l}
 \mbox{size-30 {\DFT}, breadth-first:}\\
   \quad\left \{ 
     \begin{array}{l}
      \left \{
         \begin{array}{l}
           \mbox{loop 3} \\
           \quad\mbox{size-5 direct codelet, vector size 2} \\
         \end{array}
      \right . \\
      \left \{
         \begin{array}{l}
           \mbox{loop 3} \\
           \quad\mbox{size-2 twiddle codelet, vector size 5} \\
         \end{array}
      \right . \\
    \end{array}
   \right .\\
   \quad\mbox{size-3 twiddle codelet, vector size 10}
 \end{array}
\end{array}
\]
\caption{Two possible decompositions for a size-30 DFT, both for the
arbitrary choice of DIT radices 3 then 2 then 5, and prime-size
codelets.  Items grouped by a ``$\{$'' result from the plan for a
single sub-problem.  In the depth-first case, the vector rank was
reduced to 0 as per \secref{struct:vrankr} before decomposing
sub-problems, and vice-versa in the breadth-first case.}
\label{fig:x-first}
\end{figure}

As discussed previously in \secreftwo{fftintro}{cache:ideal}, the same
Cooley--Tukey decomposition can be executed in either traditional
breadth-first order or in recursive depth-first order, where the
latter has some theoretical cache advantages.  {\FFTW} is explicitly
recursive, and thus it can naturally employ a depth-first order.
Because its sub-problems contain a vector loop that can be executed in
a variety of orders, however, {\FFTW} can also employ breadth-first
traversal.  In particular, a {\onedim} algorithm resembling the
traditional breadth-first Cooley--Tukey would result from applying
\secref{struct:rank1:Cooley--Tukey} to completely factorize the problem
size before applying the loop rule (\secref{struct:vrankr}) to reduce
the vector ranks, whereas depth-first traversal would result from
applying the loop rule before factorizing each subtransform.  These
two possibilities are illustrated by an example in \figref{x-first}.

Another example of the effect of loop reordering is a style of plan
that we sometimes call \define{vector recursion} (unrelated to
``vector-radix'' {\FFTs}~\cite{DuhamelVe90}). The basic idea is that,
if one has a loop (vector-rank $1$) of transforms, where the vector
stride is smaller than the transform size, it is advantageous to push
the loop towards the leaves of the transform decomposition, while
otherwise maintaining recursive depth-first ordering, rather than
looping ``outside'' the transform; i.e., apply the usual {\FFT} to
``vectors'' rather than numbers.  Limited forms of this idea have
appeared for computing multiple {\FFTs} on vector processors (where
the loop in question maps directly to a hardware
vector)~\cite{Swarztrauber82}.  For example, Cooley--Tukey
produces a unit \emph{input}-stride vector loop at the top-level
{\DIT} decomposition, but with a large \emph{output} stride; this
difference in strides makes it non-obvious whether vector recursion is
advantageous for the sub-problem, but for large transforms we often
observe the planner to choose this possibility.

In-place {\onedim} transforms (with no separate bit reversal pass) can
be obtained as follows by a combination {\DIT} and {\DIF} plans
(\secref{struct:rank1:Cooley--Tukey}) with transposes
(\secref{struct:rank0}).  First, the transform is decomposed via a
radix-$p$ {\DIT} plan into a vector of $p$ transforms of size $qm$,
then these are decomposed in turn by a radix-$q$ {\DIF} plan into a
vector (rank $2$) of $p \times q$ transforms of size $m$.  These
transforms of size $m$ have input and output at different
places/strides in the original array, and so cannot be solved
independently.  Instead, an indirect plan (\secref{struct:indirect})
is used to express the sub-problem as $pq$ in-place transforms of size
$m$, followed or preceded by an $m \times p \times q$ rank-$0$
transform.  The latter sub-problem is easily seen to be $m$ in-place
$p \times q$ transposes (ideally square, i.e.~$p = q$).  Related
strategies for in-place transforms based on small transposes were
described in~\cite{JohnsonBu84,Temperton91,Qian94,Hegland94};
alternating {\DIT}/{\DIF}, without concern for in-place operation, was
also considered in~\cite{Nakayama88,Saidi94}.

\subsection{The {\FFTW} planner}

Given a problem and a set of possible plans, the basic principle
behind the {\FFTW} planner is straightforward: construct a plan for
each applicable algorithmic step, time the execution of these plans,
and select the fastest one.  Each algorithmic step may break the
problem into subproblems, and the fastest plan for each subproblem is
constructed in the same way.  These timing measurements can either be
performed at runtime, or alternatively the plans for a given set of
sizes can be precomputed and loaded at a later time.

A direct implementation of this approach, however, faces an
exponential explosion of the number of possible plans, and hence of
the planning time, as $n$ increases.  In order to reduce the planning
time to a manageable level, we employ several heuristics to reduce the
space of possible plans that must be compared.  The most important of
these heuristics is \define{dynamic programming}~\cite[chapter
16]{CormenLeRi90}: it optimizes each sub-problem locally,
independently of the larger context (so that the ``best'' plan for a
given sub-problem is re-used whenever that sub-problem is
encountered).  Dynamic programming is not guaranteed to find the
fastest plan, because the performance of plans is context-dependent on
real machines (e.g., the contents of the cache depend on the preceding
computations); however, this approximation works reasonably well in
practice and greatly reduces the planning time.  Other approximations,
such as restrictions on the types of loop-reorderings that are
considered (\secref{struct:vrankr}), are described in \citeasnoun{FFTW05}.

Alternatively, there is an \define{estimate mode} that performs no
timing measurements whatsoever, but instead minimizes a heuristic cost
function.  This can reduce the planner time by several orders of
magnitude, but with a significant penalty observed in plan efficiency;
e.g., a penalty of 20\% is typical for moderate $n \lesssim 2^{13}$,
whereas a factor of 2--3 can be suffered for large $n \gtrsim
2^{16}$~\cite{FFTW05}.  Coming up with a better heuristic plan is an
interesting open research question; one difficulty is that, because
{\FFT} algorithms depend on factorization, knowing a good plan for $n$
does not immediately help one find a good plan for nearby $n$.

\section{Generating Small {\FFT} Kernels}
\label{sec:genfft}

The base cases of {\FFTW}'s recursive plans are its ``codelets,'' and
these form a critical component of {\FFTW}'s performance.  They
consist of long blocks of highly optimized, straight-line code,
implementing many special cases of the {\DFT} that give the planner a
large space of plans in which to optimize.  Not only was it
impractical to write numerous codelets by hand, but we also needed to
rewrite them many times in order to explore different algorithms and
optimizations.  Thus, we designed a special-purpose ``{\FFT}
compiler'' called \define{\genfft} that produces the codelets
automatically from an abstract description.  {\genfft} is summarized
in this section and described in more detail by~\cite{Frigo99}.

A typical codelet in {\FFTW} computes a {\DFT} of a small, fixed size
$n$ (usually, $n \leq 64$), possibly with the input or output
multiplied by twiddle factors (\secref{struct:rank1:Cooley--Tukey}).
Several other kinds of codelets can be produced by {\genfft}, but we
will focus here on this common case.

In principle, all codelets implement some combination of the
Cooley--Tukey algorithm from \eqref{ct} and/or some other {\DFT}
algorithm expressed by a similarly compact formula.  However, a high
performance implementation of the {\DFT} must address many more
concerns than \eqref{ct} alone suggests.  For example, \eqref{ct}
contains multiplications by~$1$ that are more efficient to
omit. \eqref{ct} entails a run-time factorization of~$n$, which can be
precomputed if~$n$ is known in advance. \eqref{ct} operates on complex
numbers, but breaking the complex-number abstraction into real and
imaginary components turns out to expose certain non-obvious
optimizations.  Additionally, to exploit the long pipelines in current
processors, the recursion implicit in \eqref{ct} should be unrolled and
re-ordered to a significant degree.  Many further optimizations are
possible if the complex input is known in advance to be purely real
(or imaginary).  Our design goal for {\genfft} was to keep the
expression of the {\DFT} algorithm independent of such
concerns.  This separation allowed us to experiment with various
{\DFT} algorithms and implementation strategies independently and
without (much) tedious rewriting.

{\genfft} is structured as a compiler whose input consists of the kind
and size of the desired codelet, and whose output is C~code.
{\genfft} operates in four phases: creation, simplification,
scheduling, and unparsing.

In the \define{creation} phase, {\genfft} produces a representation of
the codelet in the form of a directed acyclic graph (dag).  The dag is
produced according to well-known {\DFT} algorithms: Cooley--Tukey
(\eqref{ct}), prime-factor~\cite[page 619]{OppenheimSha89},
split-radix~\cite{Yavne68,Duhamel84,Vetterli84,Martens84,DuhamelVe90},
and Rader~\cite{Rader68}.  Each algorithm is expressed in a
straightforward math-like notation, using complex numbers, with no
attempt at optimization.  Unlike a normal {\FFT} implementation,
however, the algorithms here are evaluated symbolically and the
resulting symbolic expression is represented as a dag, and in
particular it can be viewed as a \define{linear
network}~\cite{CrochiereOp75} (in which the edges represent
multiplication by constants and the vertices represent additions of
the incoming edges).

In the \define{simplification} phase, {\genfft} applies local
rewriting rules to each node of the dag in order to simplify it.  This
phase performs algebraic transformations (such as eliminating
multiplications by~$1$) and common-subexpression elimination.
Although such transformations can be performed by a conventional
compiler to some degree, they can be carried out here to a greater
extent because {\genfft} can exploit the specific problem domain.  For
example, two equivalent subexpressions can always be detected, even if
the subexpressions are written in algebraically different forms,
because all subexpressions compute linear functions.  Also, {\genfft}
can exploit the property that \define{network transposition}
(reversing the direction of every edge) computes the transposed linear
operation~\cite{CrochiereOp75}, in order to transpose the network,
simplify, and then transpose back---this turns out to expose
additional common subexpressions~\cite{Frigo99}.  In total, these
simplifications are sufficiently powerful to derive {\DFT} algorithms
specialized for real and/or symmetric data automatically from the
complex algorithms.  For example, it is known that when the input of a
{\DFT} is real (and the output is hence conjugate-symmetric), one can
save a little over a factor of two in arithmetic cost by specializing
{\FFT} algorithms for this case---with {\genfft}, this specialization
can be done entirely automatically, pruning the redundant operations
from the dag, to match the lowest known operation count for a
real-input {\FFT} starting only from the complex-data
algorithm~\cite{Frigo99,Johnson07}.  We take advantage of this
property to help us implement real-data {\DFTs}~\cite{Frigo99,FFTW05},
to exploit machine-specific ``SIMD'' instructions
(\secref{genfft:simd})~\cite{FFTW05}, and to generate codelets for the
discrete cosine (\DCT) and sine (\DST) transforms~\cite{Frigo99,Johnson07}.
Furthermore, by experimentation we have discovered additional
simplifications that improve the speed of the generated code.  One
interesting example is the elimination of negative constants~\cite{Frigo99}:
multiplicative constants in {\FFT} algorithms often come in
positive/negative pairs, but every C~compiler we are aware of will
generate separate load instructions for positive and negative versions
of the same constants.\footnote{Floating-point constants must be
stored explicitly in memory; they cannot be embedded directly into the
CPU instructions like integer ``immediate'' constants.}  We thus
obtained a 10--15\% speedup by making all constants positive, which
involves propagating minus signs to change additions into subtractions
or vice versa elsewhere in the dag (a daunting task if it had to be
done manually for tens of thousands of lines of code).

In the \define{scheduling} phase, \genfft{} produces a topological
sort of the dag (a ``schedule'').  The goal of this phase is to find a
schedule such that a C~compiler can subsequently perform a good
register allocation.  The scheduling algorithm used by {\genfft}
offers certain theoretical guarantees because it has its foundations
in the theory of cache-oblivious algorithms~\cite{FrigoLe99} (here,
the registers are viewed as a form of cache), as described in
\secref{cache:fftw}.  As a practical matter, one consequence of this
scheduler is that {\FFTW}'s machine-independent codelets are no slower
than machine-specific codelets generated by {\SPIRAL}
\cite[Figure~3]{XiongPa01}.

In the stock {\genfft} implementation, the schedule is finally
unparsed to C.  A variation from~\cite{FranchettiKra05} implements the
rest of a compiler back end and outputs assembly code.

\subsection{SIMD instructions}
\label{sec:genfft:simd}

Unfortunately, it is impossible to attain nearly peak performance on
current popular processors while using only portable C code.  Instead,
a significant portion of the available computing power can only be
accessed by using specialized SIMD (single-instruction multiple data)
instructions, which perform the same operation in parallel on a data
vector.  For example, all modern ``x86'' processors can execute
arithmetic instructions on ``vectors'' of four single-precision values
(SSE instructions) or two double-precision values (SSE2 instructions)
at a time, assuming that the operands are arranged consecutively in
memory and satisfy a 16-byte alignment constraint.  Fortunately,
because nearly all of {\FFTW}'s low-level code is produced by
{\genfft}, machine-specific instructions could be exploited by
modifying the generator---the improvements are then automatically
propagated to all of {\FFTW}'s codelets, and in particular are not
limited to a small set of sizes such as powers of two.

SIMD instructions are superficially similar to ``vector processors'',
which are designed to perform the same operation in parallel on an all
elements of a data array (a ``vector'').  The performance of
``traditional'' vector processors was best for long vectors that are
stored in contiguous memory locations, and special algorithms were
developed to implement the {\DFT} efficiently on this kind of
hardware~\cite{Swarztrauber82,Hegland94}.  Unlike in vector
processors, however, the SIMD vector length is small and fixed
(usually 2~or~4).  Because microprocessors depend on caches for
performance, one cannot naively use SIMD instructions to simulate a
long-vector algorithm: while on vector machines long vectors generally
yield better performance, the performance of a microprocessor drops as
soon as the data vectors exceed the capacity of the cache.
Consequently, SIMD instructions are better seen as a restricted form
of instruction-level parallelism than as a degenerate flavor of vector
parallelism, and different {\DFT} algorithms are required.

The technique used to exploit SIMD instructions in {\genfft} is most
easily understood for vectors of length two (e.g., SSE2).  In this case, we view a \emph{complex} {\DFT} as a pair of \emph{real} {\DFTs}:
\begin{equation}
 \hbox{\DFT}(A + i\cdot B) =  \hbox{\DFT}(A) + i\cdot\hbox{\DFT}(B) \ ,
 \label{eq:simd-dft}
\end{equation}
where $A$ and $B$ are two real arrays.  Our algorithm computes the two
real {\DFTs} in parallel using SIMD instructions, and then it combines
the two outputs according to \eqref{simd-dft}.  This SIMD algorithm
has two important properties.  First, if the data is stored as an
array of complex numbers, as opposed to two separate real and
imaginary arrays, the SIMD loads and stores always operate on
correctly-aligned contiguous locations, even if the the complex
numbers themselves have a non-unit stride. Second, because the
algorithm finds two-way parallelism in the real and imaginary parts of
a single {\DFT} (as opposed to performing two {\DFTs} in parallel), we
can completely parallelize {\DFTs} of any size, not just even sizes or
powers of~2.

\section{Numerical Accuracy in {\FFTs}}
\label{sec:accuracy}

An important consideration in the implementation of any practical
numerical algorithm is numerical accuracy: how quickly do
floating-point roundoff errors accumulate in the course of the
computation?  Fortunately, {\FFT} algorithms for the most part have
remarkably good accuracy characteristics.  In particular, for a {\DFT}
of length $n$ computed by a Cooley--Tukey algorithm, the
\emph{worst-case} error growth is $O(\log n)$~\cite{GenSan66,Tasche00}
and the mean error growth for random inputs is only $O(\sqrt{\log
n})$~\cite{Schatzman96,Tasche00}.  This is so good that, in practical
applications, a properly implemented {\FFT} will rarely be a
significant contributor to the numerical error.

However, these encouraging error-growth rates \emph{only} apply if the
trigonometric ``twiddle'' factors in the {\FFT} algorithm are computed
very accurately.  Many {\FFT} implementations, including
{\FFTW} and common manufacturer-optimized libraries, therefore use
precomputed tables of twiddle factors calculated by means of standard
library functions (which compute trigonometric constants to roughly
machine precision).  The other common method to compute twiddle
factors is to use a trigonometric recurrence formula---this saves
memory (and cache), but almost all recurrences have errors that grow
as $O(\sqrt{n})$, $O(n)$, or even $O(n^2)$~\cite{Tasche02}, which lead to
corresponding errors in the {\FFT}.  For example, one simple
recurrence is $e^{\imagunit(k+1)\theta} = e^{\imagunit k\theta} e^{\imagunit\theta}$,
multiplying repeatedly by $e^{\imagunit\theta}$ to obtain a sequence of
equally spaced angles, but the errors when using this process grow as $O(n)$~\cite{Tasche02}.  A common improved recurrence is $e^{\imagunit(k+1)\theta} =
e^{\imagunit k\theta} + e^{\imagunit k\theta}(e^{\imagunit\theta}-1)$, where the small
quantity\footnote{In an {\FFT}, the twiddle factors are powers of
$\rootunity_n$, so $\theta$ is a small angle proportional to $1/n$ and
$e^{\imagunit\theta}$ is close to~1.}  $e^{\imagunit\theta}-1 =
\cos(\theta)-1+i\sin(\theta)$ is computed using $\cos(\theta)-1 =
-2\sin^2(\theta/2)$~\cite{Singleton67}; unfortunately, the error using
this method still grows as $O(\sqrt{n})$~\cite{Tasche02}, far worse than
logarithmic.

There are, in fact, trigonometric recurrences with the same
logarithmic error growth as the {\FFT}, but these seem more difficult
to implement efficiently; they require that a table of $\Theta(\log
n)$ values be stored and updated as the recurrence
progresses~\cite{Buneman87,Tasche02}.  Instead, in order to gain at
least some of the benefits of a trigonometric recurrence (reduced
memory pressure at the expense of more arithmetic), {\FFTW} includes
several ways to compute a much smaller twiddle table, from which the
desired entries can be computed accurately on the fly using a bounded
number (usually $< 3$) of complex multiplications. For example,
instead of a twiddle table with $n$ entries $\rootunity_n^k$, {\FFTW} can
use two tables with $\Theta(\sqrt{n})$ entries each, so that
$\rootunity_n^k$ is computed by multiplying an entry in one table (indexed
with the low-order bits of $k$) by an entry in the other table
(indexed with the high-order bits of $k$).

There are a few non-Cooley--Tukey algorithms that are known to have
worse error characteristics, such as the ``real-factor''
algorithm~\cite{Rader76,DuhamelVe90}, but these are rarely used in
practice (and are not used at all in {\FFTW}).  On the other hand,
some commonly used algorithms for type-I and type-IV discrete cosine
transforms~\cite{Swarztrauber82,PressFlaTeu92,ChanHo90} have errors
that we observed to grow as $\sqrt{n}$ even for accurate trigonometric
constants (although we are not aware of any theoretical error analysis
of these algorithms), and thus we were forced to use alternative
algorithms~\cite{FFTW05}.

To measure the accuracy of {\FFTW}, we compare against a slow {\FFT}
implemented in arbitrary-precision arithmetic, while to verify the
correctness we have found the $O(n \log n)$ self-test algorithm of
\citeasnoun{Ergun95} very useful.

\section{Generality and {\FFT} Implementations}
\label{sec:generality}

One of the key factors in {\FFTW}'s success seems to have been its
flexibility in addition to its performance.  In fact, {\FFTW} is probably
the most flexible {\DFT} library available:
\begin{itemize}
\item {\FFTW} is written in portable~C and runs well on many
  architectures and operating systems.
\item {\FFTW} computes {\DFTs} in $O(n \log n)$ time for any
  length~$n$.  (Most other {\DFT} implementations are either
  restricted to a subset of sizes or they become $\Theta(n^2)$ for
  certain values of~$n$, for example when $n$~is prime.)
\item {\FFTW} imposes no restrictions on the rank (dimensionality) of
  multi-dimensional transforms.  (Most other implementations are
  limited to one-dimensional, or at most two- and three-dimensional
  data.)
\item {\FFTW} supports multiple and/or strided {\DFTs}; for example,
  to transform a 3-component vector field or a portion of a
  multi-dimensional array.  (Most implementations support only a
  single {\DFT} of contiguous data.)
\item {\FFTW} supports {\DFTs} of real data, as well as of real
  symmetric/anti-symmetric data (also called discrete cosine/sine
  transforms).
\end{itemize}
Our design philosophy has been to first define the most general
reasonable functionality, and then to obtain the highest possible
performance without sacrificing this generality.  In this section, we
offer a few thoughts about why such flexibility has proved important,
and how it came about that {\FFTW} was designed in this way.

{\FFTW}'s generality is partly a consequence of the fact the {\FFTW}
project was started in response to the needs of a real application for
one of the authors (a spectral solver for Maxwell's
equations~\cite{Johnson2001:mpb}), which from the beginning had to run
on heterogeneous hardware.  Our initial application required
multi-dimensional {\DFTs} of three-component vector fields (magnetic
fields in electromagnetism), and so right away this meant: (i)
multi-dimensional {\FFTs}; (ii) user-accessible loops of {\FFTs} of
discontiguous data; (iii) efficient support for non-power-of-two sizes
(the factor of eight difference between $n \times n \times n$ and
$2n\times2n\times2n$ was too much to tolerate); and (iv) saving a
factor of two for the common real-input case was desirable.  That is,
the initial requirements already encompassed most of the features
above, and nothing about this application is particularly unusual.

Even for one-dimensional {\DFTs}, there is a common misperception that
one should always choose power-of-two sizes if one cares about
efficiency.  Thanks to {\FFTW}'s codelet generator, we could afford to
devote equal optimization effort to any $n$ with small factors (2, 3,
5, and 7 are good), instead of mostly optimizing powers of two like
many high-performance {\FFTs}.  As a result, to pick a typical example
on the 3~GHz Core Duo processor of \figref{fftwnr}, $n=3600=2^4 \cdot
3^2 \cdot 5^2$ and $n=3840 = 2^8 \cdot 3 \cdot 5$ both execute faster
than $n=4096=2^{12}$.  (And if there are factors one particularly
cares about, one can generate codelets for them too.)

One initially missing feature was efficient support for large prime
sizes; the conventional wisdom was that large-prime algorithms were
mainly of academic interest, since in real applications (including
ours) one has enough freedom to choose a highly composite transform
size.  However, the prime-size algorithms are fascinating, so we
implemented Rader's $O(n \log n)$ prime-$n$ algorithm~\cite{Rader68} purely for fun,
including it in {\FFTW}~2.0 as a bonus feature.  The response was
astonishingly positive---even though users are (probably) never
\emph{forced} by their application to compute a prime-size {\DFT}, it
is rather inconvenient to always worry that collecting an unlucky
number of data points will slow down one's analysis by a factor of a
million.  The prime-size algorithms are certainly slower than
algorithms for nearby composite sizes, but in interactive
data-analysis situations the difference between 1~ms and 10~ms means little,
while educating users to avoid large prime factors is hard.

Another form of flexibility that deserves comment has to do with a
purely technical aspect of computer software.  {\FFTW}'s
implementation involves some unusual language choices internally
({\genfft} is written in Objective Caml, a functional language
especially suited for compiler-like programs), but its user-callable
interface is purely in C with lowest-common-denominator datatypes
(arrays of floating-point values).  The advantage of this is that
{\FFTW} can be (and has been) called from almost any other programming
language, from Java to Perl to Fortran~77.  Similar
lowest-common-denominator interfaces are apparent in many other
popular numerical libraries, such as LAPACK~\cite{LAPACK}.  Language
preferences arouse strong feelings, but this technical constraint means
that modern programming dialects are best hidden from view for a
numerical library.

Ultimately, very few scientific-computing applications should have
performance as their top priority.  Flexibility is far more important,
because one wants to be limited only by one's imagination, rather than
by one's software, in the kinds of problems that can be studied.

\section{Concluding Remarks}

It is unlikely that many readers of this chapter will ever have to
implement their own fast Fourier transform software, except as a
learning exercise.  The computation of the {\DFT}, much like basic
linear algebra or integration of ordinary differential equations, is
so central to numerical computing and so well-established that robust,
flexible, highly optimized libraries are widely available, for the
most part as free/open-source software.  And yet there are many other
problems for which the algorithms are not so finalized, or for which
algorithms are published but the implementations are unavailable or of
poor quality.  Whatever new problems one comes across, there
is a good chance that the chasm between theory and efficient
implementation will be just as large as it is for {\FFTs}, unless
computers become much simpler in the future.  For readers who
encounter such a problem, we hope that these lessons from {\FFTW} will be
useful:

\begin{itemize}
\item Generality and portability should almost always come first.
\item The number of operations, up to a constant factor, is less important than the order of operations.
\item Recursive algorithms with large base cases make optimization easier.
\item Optimization, like any tedious task, is best automated.
\item Code generation reconciles high-level programming with low-level performance.
\end{itemize}

We should also mention one final lesson that we haven't discussed in
this chapter: you can't optimize in a vacuum (or you end up
congratulating yourself for making a slow program slightly faster).
We started the {\FFTW} project after downloading a dozen {\FFT}
implementations, benchmarking them on a few machines, and noting how
the winners varied between machines and between transform sizes.
Throughout {\FFTW}'s development, we continued to benefit from
repeated benchmarks against the dozens of high-quality {\FFT} programs
available online, without which we would have thought {\FFTW} was
``complete'' long ago.

\section*{Acknowledgements}

SGJ was supported in part by the Materials Research Science and
Engineering Center program of the National Science Foundation under
award DMR-9400334; MF was supported in part by the Defense Advanced
Research Projects Agency (DARPA) under contract No. NBCH30390004.  We
are also grateful to Sidney Burrus for the opportunity to contribute
this chapter, and for his continual encouragement---dating back to his
first kind words in 1997 for the initial {\FFT} efforts of two
graduate students venturing outside their fields.

\bibliography{fft}
\bibliographystyle{ieeetr}

\end{document}